\documentstyle{article}
\begin{document}

\newtheorem{satz}{Theorem}[section]
\newtheorem{lem}[satz]{Lemma}
\newtheorem{prop}[satz]{Proposition}
\newtheorem{kor}[satz]{Corollary}
\noindent
ON SUPERHEIGHT CONDITIONS FOR THE AFFINENESS OF OPEN SUBSETS
\par
\bigskip
\noindent
\centerline{Holger Brenner, Bochum}
\par
\bigskip
\noindent
\begin{abstract}
\noindent
In this paper we consider the open complement $U$ of a hypersurface
$Y=V({\bf a})$ in an affine scheme $X$.
We study the relations between the
affineness of $U$, the intersection of $Y$ with closed subschemes,
the property that every closed surface in $U$ is affine,
the property that every analytic closed surface is Stein
and the superheight of a defining ideal $\bf a$.
\end{abstract}
\par\bigskip
\noindent
\section{Introduction}
\par
\bigskip
\noindent
Let $A$ be a noetherian ring with an ideal ${\bf a} \subseteq A$
and let $X={\rm Spec}\, A$, $Y=V({\bf a})$. We consider the complement
$U=D({\bf a})=X-Y$. The purpose of this paper is to find geometric
conditions for $U$ to be affine.
It is well known that if $U$ is affine then $Y$ must be a
hypersurface, i.e. ${\rm ht}\, {\bf a} \leq 1$, see prop. 2.4.
Note that the converse is by no means true, yet the height condition
on $Y$ has a stronger generalization based on the following simple
observation.
\par\smallskip\noindent
Let $X'={\rm Spec}\, A'$ be another affine scheme and
$f:X' \longrightarrow X$ be a morphism corresponding to the
ring homomorphism $A \longrightarrow A'$. If $U=D({\bf a}) \subseteq X$
is affine then the preimage $U'=f^{-1}(U)=D({\bf a}A')$
is also an affine scheme.
Since $Y'=f^{-1}(Y)=V({\bf a}A')$ this means that the height condition must
also hold for all minimal primes of the extended ideal ${\bf a}'={\bf a}A'$.
\par\smallskip\noindent
It is a general observation, first studied by Neeman (\cite{neeman}),
that the non\--affineness can often be shown
by giving a ring homomorphism violating the height condition on
the extended ideal of $\bf a$. In order to illustrate this technique
we give the following example.
\par\bigskip
\noindent
{\bf Example 1.1.}
Let $K$ be a field, $A=K[R,S,T,Z]/(RS-TZ)$ and $X={\rm Spec}\, A$.
$A$ is a normal three\--dimensional domain, the ideal
${\bf a}=( R ,T ) $ is prime of
height one. Let $Y=V(R,T)$.
Under the reduction $A \longrightarrow A/(Z,S)=K[R,T]$
the extended ideal is $(R,T)$ in $K[R,T]$ which is of height two.
Since the complement of a point in the plane is not affine it follows
that $U=X-Y$ cannot be affine either.
\par\bigskip
\noindent
In this paper we study the connection between the affineness of $U=D({\bf a})$
and the property that the codimension of $Y'$ under every ring homomorphism
is restricted by one. This property can be expressed in terms
of the superheight of the ideal, namely ${\rm supht}\, {\bf a} \leq 1$.
This notion was first introduced by Hochster in 1975 \cite{hochsterbig}.
We give an intrinsic definition of superheight depending only on the open set
$D({\bf a})$ (not on the ideal) so that the notion of superheight
can be extended to arbitrary schemes (2).
\par\smallskip\noindent
In section (3) we describe situations where the affineness can be obtained
from superheight conditions. We show that the affineness of $D({\bf a})$
is equivalent
to the property that for all ring homomorphisms $A \longrightarrow A'$
where $A'$ is a Krull domain the height of the extended ideal is $\leq 1$.
For a noetherian domain we characterize the affineness in terms of
finite superheight under the additional condition that the ring of
global sections is finitely generated, generalizing a result of
Neeman \cite{neeman}. Furthermore, in the two dimensional case and in the case of
monoid rings the affineness can be read off directly from the behaviour
of the height in only one special ring extension.
\par\smallskip\noindent
In section (4) we consider schemes of finite type over the complex numbers
$\bf C$ and define the notion of analytic superheight and compare it with the
algebraic notions of superheight.
It will turn out that if $U$ is Stein as a complex space then
the analytic and the algebraic superheight is one.
We recover the result of Bingener and Storch (\cite{binsto1})
that, under the condition
that the ring of global sections is finitely generated, affineness and Stein
are the same.
\par\smallskip\noindent
In section (5) we consider finitely generated $K-$Algebras and relate
the superheight one condition to the property that
every closed subscheme of $U$ of dimension $ \leq 2$ is affine.
We show that in the complex case this property is equivalent to
the property that any closed analytic surface in $U$ is Stein.
The question whether this last property implies the Stein property
for $U$ is the so\-- called hypersurface (or hypersection) problem
answered negatively by Coltoiu and Diederich
(\cite{diederich1},\cite{diederich2}).
\par\smallskip\noindent
Finally, in section (6), we give two classes of examples
of non\--affine open subsets with superheight one.
The first class is constructed from certain curves on
smooth projective surfaces,
using the intrinsic characterization of superheight
and a criterion
\`{a} la ``Riemannscher Fortsetzbarkeitssatz'' for superheight one.
The other is built from non-torsion divisor classes of a local
two dimensional normal ring, related to a construction of Rees and
yielding counterexamples to the hypersurface problem.
\par
\bigskip
\noindent
\section{The superheight of an ideal and of a scheme}
\par
\bigskip
\noindent
Let ${\bf a} \subseteq A$ be an ideal in a commutative ring and
$A \longrightarrow A'$ a ring homomorphism.
The extended ideal
${\bf a}A'$ describes the preimage of the
open set $D({\bf a})$ under the mapping
${\rm Spec}\, A' \longrightarrow {\rm Spec}\, A$.
The height of an ideal ${\bf a} \subset A$ is defined as the minimal
height of a minimal prime of $\bf a$. The maximal height of
the minimal primes is called the big height or the altitude of
$\bf a$.
We put ${\rm ht}\, ({\bf A}) = 1$ in case $A \neq 0$, otherwise
$=0$.
\par\bigskip
\noindent
{\bf Definition.}
For an ideal ${\bf a} \subseteq A$ in a commutative ring we
call
\par\medskip\noindent
${\rm supht}\, {\bf a}=
{\rm max}\, \{   {\rm ht}\, {\bf a}A':\, A \longrightarrow A' \mbox{ with }
A' \,\, {\rm noetherian} \}$
the {\it superheight} of $\bf a$ or the noetherian superheight.
\par\medskip\noindent
${\rm supht}^{\rm fin}\, {\bf a}=
{\rm max}\, \{   {\rm ht}\, {\bf a}A': \, A \longrightarrow A' \mbox{ with }
A'\,\, {\rm of\, finite\, type} \}$
the {\it finite superheight}
of $\bf a$.
\par\medskip\noindent
${\rm supht}^{\rm krull}\, {\bf a}
={\rm max}\, \{   {\rm ht}\, {\bf a}A': \, A \longrightarrow A' \mbox{ with }
A'\,\, {\rm Krull \, domain} \}$
the superheight with respect to Krull domains.
\par\bigskip
\noindent
This notion goes back to Hochster and was developed in connection with the
direct summand conjecture, \cite{hochsterbig}.
This conjecture states that a
local regular ring $A$ is a direct summand in any finite extension
$A \subseteq B$. 
The conjecture is known to be true if $A$ contains a field.
In general it is equivalent to the monomial conjecture, which can be stated
as a proposition about the superheight of an ideal, namely that the
ideal $(X_1,...,X_n)$ in
$${\bf Z}[X_1,...,X_n,Y_1,...,Y_n]/
((X_1 \cdot ... \cdot X_n)^{k} -Y_1X_1^{k+1}-...-Y_nX_n^{k+1})$$
has superheight $n-1$ (for every $k \in {\bf N})$,
see \cite{hochsterdim}, \cite{hochstercan} and below for the treatment
of the two dimensional case ($n=2)$ via affineness.
There are some important results of Koh (\cite{koh1}, \cite{koh2})
on superheight which we will use in the following.
\begin{prop}\hspace{-0.5em}{\bf .}
Let $\bf a$ be an ideal in a noetherian ring $A$.
Then the following statements are true
{\rm (}in {\rm (2)} and {\rm (3)} suppose ${\bf a } \neq A$ {\rm )}.
\par\smallskip\noindent
{\rm (1)}
$$ {\rm ht}\, {\bf a} \leq {\rm bight}\, {\bf a}
\leq {\rm supht}^{\rm fin}\, {\bf a}
\leq {\rm supht}\, {\bf a}
\leq {\rm ara}\, {\bf a} \, .$$
\par\smallskip\noindent
{\rm (2)}
The finite superheight equals
\begin{eqnarray*}
{\rm supht}^{\rm fin}\, {\bf a} &= & {\rm max}\, \{ {\rm ht}\, {\bf m}:\, 
{\bf m} \mbox{ is a maximal ideal of } A',\,
A \longrightarrow A' \mbox{ is of finite } \\
& & \ \ \ \ \ \ \ \ \ \ \ \ \ \ \ \ \ \ \ \ \ \ \ \ \ \ \ \ \ \ \ \ \ \ \ \ \
\ \ \ \ \ \mbox{ type and }
V({\bf a}A')=V({\bf m}) \}
\end{eqnarray*}
\par\smallskip\noindent
{\rm (3)}
The superheight equals
\begin{eqnarray*}
{\rm supht}\, {\bf a} &= & {\rm max}\, \{ {\rm dim}\, A':\, 
A' \mbox{ is a noetherian local complete normal domain, } \\
& & \ \ \ \ \ \ \ \ \ \ A \longrightarrow A'
\mbox{ is a ring homomorphism with } V({\bf a}A')=V({\bf m}) \}
\end{eqnarray*}
\end{prop}
{\it Proof}.
(1).
The first and the third inequality are clear, the second is proved below.
${\rm ara}\, {\bf a}$ denotes the minimal number of functions $f_1,...,f_k$
with $V({\bf a})=V(f_1,...,f_k)$, so the forth inequality
follows from the general Krull Hauptidealsatz,  \cite{eisenbud}, Theorem 10.2.
\par\smallskip\noindent
(2) and (3). Let ${\bf a}A' \subseteq {\bf p}$ be a minimal prime ideal in $A'$,
and ${\bf p}_1,...,{\bf p}_r$ the others.
Prime avoidance shows that there exists
$f \not\in {\bf p}$ and $f \in {\bf p}_i$ for $ i=1,...,r$.
After the change $A' \longrightarrow A'_f$ the prime ideal ${\bf p}A'_f$
is the only minimal prime over ${\bf a}A'_f$, and the height is the 
height of $\bf p$. This proves the second inequality of (1).
If the height of $\bf p$ is taken over the prime $\bf q$ of height zero,
we get modulo $\bf q$ a domain.
So in both cases we can restrict to morphisms where $A'$ is a domain and
$V({\bf a}A')=V({\bf p})$ irreducible.
\par\smallskip\noindent
We show that for a prime ideal $\bf p$ of height $n$
in a noetherian ring there exists a residue class domain where
$\bf p$ extends geometrically (as a radical) to a maximal ideal of height $n$.
If $\bf p$ is maximal we are done, so let ${\bf q}$ be a direct
prime over $\bf p$.
Let $x \not\in {\bf p},\, x \in {\bf q}$.
Then ${\bf q}$ is modulo $x$ a minimal prime over ${\bf p}A'/x$
(and over ${\bf a}A'/x$), and we have in $A'/x$ the relations
(with ${\bf q}'={\bf q}(A'/x)$)
$$ {\rm ht}\, {\bf q}' =
{\rm dim}\, (A'/x)_{{\bf q}'}
={\rm dim}\, (A'_{\bf q}/x)
\geq {\rm dim}\, (A'_{\bf q}) -1
\geq {\rm ht}\, {\bf p} \, .$$
cf. \cite{eisenbud}, Cor. 10.9.
Sucessively we arrive at a maximal ideal. This proves (2).
\par\smallskip\noindent
(3). Localization at $\bf p$ yields a local ring, there
the extended ideal describes geometrically exactly the maximal
ideal, and the dimension is the superheight of $\bf a$.
Under completion the dimension does not change
and in considering a component of maximal
dimension we have a complete local domain $R$.
$R$ is excellent, and therefore \cite{EGAIV}, 7.6.2,
its normalization is again noetherian complete and local of the
same dimension. \hfill $\Box$
\begin{satz}\hspace{-0.5em}{\bf .}
Let $A$ be a noetherian ring and ${\bf a} \subseteq A$ an ideal. Then
\begin{eqnarray*}
& & {\rm supht}^{\rm fin}\, {\bf a} \\
& & ={\rm sup}\, \{ {\rm bight}\,  {\bf a}A':\,
A' \, \mbox{is the normalization of a residue class domain} \}.
\end{eqnarray*}
\end{satz}
{\it Proof}. See \cite{koh1}. \hfill $\Box$
\begin{satz}\hspace{-0.5em}{\bf .}
Let $K$ be a field and $A$ a finitely generated $K-$Algebra,
${\bf a} \subseteq A$ an ideal.
Then ${\rm supht}^{\rm fin}\, {\bf a}={\rm supht}\, {\bf a}$.
\end{satz}
{\it Proof}.
See \cite{koh2}, theorem 1.  \hfill $\Box$
\par
\bigskip
\noindent
We extend the notion of superheight to an arbitrary scheme.
\par
\bigskip
\noindent
{\bf Definition.}
Let $X$ be a scheme. The {\it superheight} of $X$ is the biggest number $d$
such that there exists
\par\smallskip\noindent
i) a noetherian affine scheme $T$ with a closed point $P \in T$ of height $d$.
\par\smallskip\noindent
ii) an affine morphism $f:T- \{P\} \longrightarrow X$.
\par\bigskip\noindent
If $X$ is a variety over a field $K$, we call the same number, under the
restriction that $T$ be an affine variety, the finite superheight of $X$.
\par\bigskip
\noindent
{\bf Remark.}
In determining the superheight of a scheme
one may only look at local complete normal
noetherian domains $T={\rm Spec}\, A$. For this, first localize at $P$ and then
do the same steps as in the proof of 2.1 (3).
\par\smallskip\noindent
If $X$ is empty we have ${\rm supht}\, X = 0$, because then
$T -\{ P\}$ has to be empty, hence ${\rm dim}\, T =0$.
On the other hand, a non-empty scheme $X$ has superheight $\geq 1$.
For a point ${\rm  Spec}\, K \longrightarrow X$ ($K$ a field)
the morphism 
$$ {\rm Spec}\, K[Y]_{(Y)} \supseteq D(Y)={\rm Spec}\, K(Y)
\longrightarrow {\rm Spec}\, K \longrightarrow X $$
is affine and $T={\rm Spec}\, K[Y]_{(Y)}$ is one\--dimensional.
\par\smallskip\noindent
If $X$ is affine, we have ${\rm supht}\, X \leq 1$, because
in this case the affineness of $T-\{ P \} \longrightarrow X$
implies the affineness of $T-\{ P\}$. The following proposition,
which is the starting point of this whole subject,
shows that ${\rm dim}\, T \leq 1$.
\begin{prop}\hspace{-0.5em}{\bf .}
Let $X$ be a noetherian separated scheme and $U \subseteq X$ an affine
open subscheme. Then every component of $Y=X-U$ has codimension $\leq 1$.
The same is true for $X={\rm Spec}\, A$ where $A$ is a Krull domain.
\end{prop}
{\it Proof}.
Let $\eta$ be the generic point of a component of $Y$, $A={\cal O}_\eta$.
Since $X$ is separated, $i:{\rm Spec}\, A \hookrightarrow X$ is an affine
morphism and thus $D(\eta)=i^{-1}(U)$ is again affine.
So we have to show that in a local noetherian ring $A$ the complement
of the closed point is affine only in case ${\rm dim}\, A \leq 1$.
We may assume that $A$ is a domain.
The normalization $A_{\rm nor}$ of $A$ is a semilocal Krull domain
(see \cite{nagloc}), so we are led to a local Krull domain $A$.
But for a Krull domain with ${\rm dim}\, A \geq 2$ we have
$\Gamma(D({\bf m}),{\cal O}_X)=A$, hence $D({\bf m})$
is not affine. \hfill $\Box$
\par
\bigskip
\noindent
The assumption in the following criterion for ${\rm supht}\,(X)\leq 1$ says
that $X$ satisfies as target the ``Riemannscher
Fortsetzbarkeitssatz''. 
\begin{lem}\hspace{-0.5em}{\bf .}
Let $X$ be a noetherian separated scheme satisfying the following property:
\par\smallskip\noindent
For any normal noetherian scheme $T$ and any closed point $P \in T$ with
${\rm ht}\, P \geq 2$, every morphism  $T- \{P\} \longrightarrow X$
is extendible to $T$.
\par\smallskip\noindent
Then ${\rm supht}\,(X)\leq 1$.
\par\smallskip\noindent
If $X$ is quasi\--affine, the converse is also true.
\end{lem}
{\it Proof}.
Let $T$ be affine.
If $X$ is separated, an affine morphism
$f: T- \{P \} \longrightarrow X$ with
${\rm ht}\, P \geq 2$ cannot be extended
to the whole of $T$.
An extension $\tilde{f}:T \longrightarrow X$ would be an
affine morphism, and for an affine open neighbourhood
$\tilde{f}(P) \in V$
the sets $\tilde{f}^{-1}(V)$ and $f^{-1}(V)$ must both be affine.
But $\tilde{f}^{-1}(V)=f^{-1}(V) \cup \{ P \}$ and $P$ is a point of height
$\geq 2$, so this is not possible.
Therefore an affine morphism $T-\{ P \} \longrightarrow X$ with ${\rm ht}\, P \geq 2$
and $T$ normal contradicts the assumption.
\par\smallskip\noindent
Let $X \subseteq {\rm Spec}\, A$ be quasi\--affine with superheight $\leq 1$
and $f: T - \{P\} \longrightarrow X$ a morphism with $T$ normal and affine,
${\rm ht} \,(P) \geq 2$. $f$ is not affine, but there is an affine extension
$\tilde{f}: T \longrightarrow {\rm Spec}\, A$, corresponding to the
ring homomorphism
$A \longrightarrow \Gamma (T-\{P\},{\cal O}_T)=\Gamma (T,{\cal O}_T)$.
If $\tilde{f}(P) \not\in X$, $f$ would be the restriction of
$\tilde{f}$ on $X$, hence affine. So $\tilde{f}(P) \in X$ and $f$ is extendible
as a mapping to $X$. \hfill $\Box$
\par\bigskip
\noindent
As the following proposition shows, the superheight of an ideal $\bf a$
and the superheight of the open set $D({\bf a})$ coincide.
\begin{prop}\hspace{-0.5em}{\bf .}
For an ideal ${\bf a} \subseteq A$ the equality
${\rm supht}\, D({\bf a})={\rm supht}\, {\bf a}$ holds.
\end{prop}
{\it Proof}.
Let $A \longrightarrow R$ be a ring homomorphism in a local
normal noetherian domain of dimension $m={\rm supht}\, {\bf a}$
with $V({\bf a}R)=V({\bf m}_R)=\{ P \}$. Let $U=D({\bf a})$.
Then the mapping
$f^{-1}(U)={\rm Spec}\, R - \{ P \} \longrightarrow U$ is affine and
therefore ${\rm supht}\, U \geq {\rm supht}\, {\bf a}$.
\par\smallskip\noindent
For the converse inequality let
$f:T- \{ P\} \longrightarrow U$ be an affine morphism where
$T$ is local normal noetherian and $d={\rm dim}\, T = {\rm supht}\, U$.
If $d=0$, there is nothing to show.
If $d=1$, it follows that $U$ is not empty.
Let ${\bf q} \in U$ be a prime ideal
and consider $A \longrightarrow A_{\bf q}=A' $.
Then we have ${\bf a}A_{\bf q}=A_{\bf q}$ and
therefore ${\rm supht}\, {\bf a} \geq 1$ follows
from the definition.
\par\smallskip\noindent
So let $d \geq 2$.
Since $T={\rm Spec}\, R$ is normal, we have
$\Gamma(T-\{ P \} ,{\cal O}_T)=\Gamma(T,{\cal O}_T)=R$
and $f$ corresponds to the global ring homomorphismus
$A \longrightarrow R$, so the mapping is extendible
to a mapping $\bar{f}: {\rm Spec}\, R \longrightarrow {\rm Spec}\, A$.
$\bar{f}(P) \in U$ is not possible, for otherwise
the mapping would be extendible as a mapping into $U$,
but this is excluded by the proof of the previous lemma.
So under $A \longrightarrow R$ the extended ideal describes
$V({\bf a}R) =\{ P\}$ and therefore
${\rm supht}\, {\bf a} \geq d={\rm supht}\, U$. \hfill $\Box$
\par\bigskip
\noindent
We gather together some properties of the superheight of a scheme.
\begin{prop}\hspace{-0.5em}{\bf .}
{\rm (1)}
For an affine morphism $X' \longrightarrow X$ we have
${\rm supht}\, X' \leq {\rm supht}\, X$.
\par\smallskip\noindent
{\rm (2)}
The superheight of $X$ equals the maximum of the superheights
of the irreducible components of $X$.
\par\smallskip\noindent
{\rm (3)}
Suppose $X$ noetherian. For $Y \subseteq X$ closed and
$U=X-Y$ we have
$$ {\rm supht}\, X \leq {\rm supht}\, Y + {\rm supht}\, U \, .$$
\par\smallskip\noindent
{\rm (4)}
If $X=U \cup V$ with $U,V$ open, we have
$ {\rm supht}\, X \leq {\rm supht}\, V + {\rm supht}\, U $.
\par\smallskip\noindent
{\rm (5)}
${\rm supht}\, X \leq {\rm dim}\, X +1$.
\par\smallskip\noindent
{\rm (6)} For a noetherian separated scheme $X$ we have
$ {\rm supht}\, X \leq  {\rm cd}\, X +1 $
{\rm (}${\rm cd}$ denotes the cohomological dimension of $X$
in the sense of R. Hartshorne,
meaning the maximal number $ n \in {\bf N}$
such that there is a quasicoherent sheaf $\cal F$ on $X$
with $H^n(X,{\cal F}) \neq 0$.{\rm )}
\end{prop}
{\it Proof}.
(1) is clear.
(2). Let $f:T-\{ P \} \longrightarrow X$ be affine with $T$ irreducible.
The image of $T$ lies in a component $X_i$ of $X$ and this
component must have the superheight of $X$.
\par\smallskip\noindent
(3). Let $f: T \supseteq T-\{ P \} \longrightarrow  X$ be an affine
morphism with $T={\rm Spec}\, R$, $R$ being a noetherian local complete
domain,
and with ${\rm dim}\, T = {\rm supht}\, X$.
Set $f^{-1}(Y)=V({\bf a})-\{ P\}$ with an ideal ${\bf a} \subseteq R$.
On one hand, we have ${\rm dim}\, V({\bf a}) \leq {\rm supht}\, Y$
as is shown by the restriction $V({\bf a})-\{ P \}=f^{-1}(Y) \longrightarrow Y$.
On the other hand, the restriction
$D({\bf a})=f^{-1}(U) \longrightarrow U$ is also affine and so
${\rm bight}\,  {\bf a} \leq {\rm supht}\, D({\bf a}) \leq {\rm supht}\, U$.
Since $R$ is complete, $R$ is catenary, (see \cite{eisenbud}, Cor. 18.10)
and so for a minimal prime $\bf p$
of $\bf a$ we have the inequalities
\begin{eqnarray*}
{\rm supht}\, X ={\rm dim}\, R &
=& {\rm dim}\, R/{\bf p}+{\rm ht}\, {\bf p} \\
& \leq & {\rm dim}\, R/{\bf a} +{\rm bight}\,  {\bf a} \\
& \leq & {\rm supht}\, Y + {\rm supht}\, (X-Y)
\end{eqnarray*}
(4). Let $X=U \cup V$. Then $Y=X-U$ is a closed subset
of $V$ leading to ${\rm supht}\, Y \leq {\rm supht}\, V$
and the statement follows from (3).
\par\smallskip\noindent
(5). We do induction on the dimension, the beginning is clear.
Because of (2) we may assume that $X$ is irreducible of dimension
$d$. For a non\--empty
open affine subset $U$, (3) yields
${\rm supht}\, X \leq {\rm supht}\, U+ {\rm supht}\, (X-U)
\leq 1+d $.
\par\smallskip\noindent
{\rm (6)}.
${\rm supht}\, X=0$ if and only if $X= \emptyset $. In this case
${\rm cd}\, X=-1$. So suppose ${\rm supht}\, X \geq 1$.
If $T$ is a local noetherian affine scheme of dimension $d \geq 1$
a theorem of Grothendieck says that $H^{d}_{{\bf m}}(T, {\cal O}) \neq 0$.
The natural map of local cohomology
$H^{i-1}(T-\{ P \}, {\cal O}) \longrightarrow H^{i}_{{\bf m}}(T, {\cal O})$
is bijective for $i \geq 2$ and surjective for $i=1$.
Thus $H^{d-1} (T-\{ P \}, {\cal O}) \neq 0$.
If $f: T-\{ P \} \longrightarrow X$ is affine and $d={\rm supht}\, X$
it follows that $H^{d-1}(X,f_\ast {\cal O}) \neq 0$ and
${\rm cd}\, X \geq d-1 = {\rm supht}\, X -1$. 
This gives also another proof of {\rm (5)} and of 2.4.  \hfill $\Box$
\par\bigskip
\noindent
{\bf Example 2.1.}
Let $Y$ be a projective variety of dimension $d$. The mapping of a
punctured affine cone $X- \{P\} \longrightarrow Y$ is affine,
hence the superheight of $Y$ is $\geq d+1$ and equality must hold
because of (5).
It is reasonable to ask whether maximal possible superheight
\---the existence of such an affine cone\---
ensures for a normal separated variety projectivity.
A result of Kleiman states that a normal separated variety is proper
if and only if the cohomological dimension is maximal,
\cite{kleiman}.
\begin{kor}\hspace{-0.5em}{\bf .}
Let $X$ be a scheme with ${\rm supht}\, X \leq d$.
Then the complement of $X$ in any open embedding $X \subseteq X'$ with
$X'$ noetherian and separated has codimension $\leq d$.
\end{kor}
Proof.
For an affine subset $U$ of $X'$ the morphism $U \cap X \hookrightarrow X$
is affine, so $U \cap X$ fulfills the assumption as well.
Since the conclusion is local, we may assume $X'$ to be affine.
Thus the statement follows from 2.6. \hfill $\Box$
\par
\bigskip
\noindent
\section{Affineness and superheight one}
\par
\bigskip
\noindent
Let ${\bf a}=(f_1,...,f_n) \subseteq A$ be an ideal in a commutative ring,
$U=D({\bf a}) \subseteq {\rm Spec}\, A=X$ and
$B=\Gamma(U,{\cal O}_X)$ the ring of global sections on $U$.
In this situation we have an open embedding
$U=D({\bf a}B) \hookrightarrow {\rm Spec}\, B$,
and $U$ is affine if and only if ${\bf a}B$ is the unit ideal.
In this case we have
$1=q_1f_1+...+q_nf_n$ with $q_i \in \Gamma(U,{\cal O}_X)$,
and the functions yields a closed embedding
$(q_1,...,q_n):U \hookrightarrow {\rm Spec}\, A[T_1,...,T_n]$, showing
by the way that in the affine case $B$ is an $A-$algebra of finite type.
\par\smallskip\noindent
If this is not the case, the height of this extended ideal is
larger than one.
\begin{satz}\hspace{-0.5em}{\bf .}
Let $A$ be a noetherian ring and $\bf a$ an ideal, $U=D({\bf a})$.
Then $U$ is affine if and only if ${\rm supht}^{\rm krull}\, {\bf a} \leq 1$.
\end{satz}
{\it Proof}.
If $U$ is affine and $A \longrightarrow A'$ is a ring homomorphism, where
$A'$ is a Krull domain, then the preimage $U'=D({\bf a}A')$ is affine
and $V(aA')$ has codimension $\leq 1$, see 2.4.
\par\smallskip\noindent
So suppose $U$ is not affine.
Since a noetherian scheme is affine if and only if all its (reduced)
components are affine (see \cite{haramp}, II.1.4) we find
$A \longrightarrow A'$, where
$A'$ is a noetherian domain and where $D({\bf a}A')$ is not affine.
So we may assume that $A$ is a domain.
Consider the normalization
$\varphi: {\rm Spec}\, A_{\rm nor} \longrightarrow {\rm Spec}\, A$.
If ${\bf a}=(f_1,...,f_n)$ and
$\varphi^{-1}(U)$ were affine, there would exist
$q_i \in \Gamma(\varphi^{-1}(U),{{\cal A}_{\rm nor}})$
with $q_1f_1+...+q_nf_n=1$.
But these functions are already
defined on the corresponding open set in a finite extension
$A \subset B \subset A_{\rm nor}$,
and the theorem of Chevalley (\cite{EGAIV} or \cite{haramp}, II.1.5
shows that $U$ itself would be affine.
\par\smallskip\noindent
So we may assume that $A$ is a Krull domain.
For an open subset $W$ in ${\rm Spec}\, A$ of a Krull domain
the ring of global sections is given by the intersection of discrete
valuation domains, 
$$\Gamma(W,{\cal O}_X)= 
\bigcap_{{\rm ht}\,({\bf p})=1,\, {\bf p} \in W} A_{\bf p}\,  .$$
From this we see that the ring of global sections $B=\Gamma(U,{\cal O}_X)$
is again a Krull domain.
We have $U=D({\bf a}) \cong D({\bf a}B) \subseteq {\rm Spec}\, B$,
and ${\bf a}B$ is not the unit ideal.
On the other hand, we have $B=\Gamma(U,{\cal B})$, and this can only
hold if $U$ contains all prime ideals of height one of the Krull domain
$B$. For if ${\bf p}$ of height one is not in $U$, let $p$ be a
generator of the maximal ideal in the discrete valuation ring
$B_{\bf p}$ and let $q=1/p$. Let ${\bf q}_1,...,{\bf q}_m$ be the other poles
of $q$. We find $f \in B$ with $f \not\in {\bf p},\, f \in {\bf q}_i $
for $i=1,...,m$. Then for all $n$ big enough the function $f^nq$ has its only
pole in $\bf p$ and is defined on $U \subseteq D({\bf p })$.
So we conclude that ${\bf a}B$ has height $\geq 2$. \hfill $\Box$
\par
\bigskip
\noindent
Under additional conditions on the ring of global sections the superheight
condition for smaller classes of rings guarantees affineness. The
following result can also be found in \cite{neeman} in the case that $A$ is
normal and of finite type over a field.
\begin{satz}\hspace{-0.5em}{\bf .}
Let $A$ be a noetherian domain and ${\bf a}$ an ideal, $U=D({\bf a})$.
Then $U$ is affine if and only if the ring of
global sections $\Gamma(U,{\cal O}_X)$ is of finite type over $A$ and
${\rm supht}^{\rm fin}\, {\bf a} \leq 1$. 
\end{satz}
{\it Proof}.
If $U$ is affine, it is known that $B=\Gamma(U,{\cal O}_X)$
is finitely generated over $A$, so suppose
$U$ is not affine with a finitely generated ring $B$ of global sections.
$B$ is a noetherian domain and the extended ideal
is not the unit ideal, but $U \cong D({\bf a}B)$
contains all prime ideals of height one of $B$.
For if ${\bf p}=(f_1,...,f_m)$ is a prime ideal in $B$ of height one
there is a function
$f \in {\bf p}$ with ${\rm Rad} (f)={\bf p}B_{\bf p}$.
This yields equations $f_i^n=(a_i/r_i)f$ with $r_i \not\in {\bf p}$.
With $r=r_1\cdot ... \cdot r_m$ we may write
$f_i^n=(b_i/r)f$ or $r/f=b_i/f_i^n$, showing that this is a function
defined on $D({\bf p})$ not belonging to $B$, since otherwise
$f (r/f) = r \in {\bf p}$.
So we have ${\rm height}\, {\bf a}B \geq 2$ and
${\rm supht}^{\rm fin}\, {\bf a} \geq 2$. \hfill $\Box$
\par
\bigskip
\noindent
In deciding whether an open subset of an affine scheme is again affine,
one may look at the ring of global section and the height
of the extended ideal in it. If this ideal is the unit ideal, $U$ is affine
and the superheight is one. If this is not the case,
the ring of global sections is just one candidate among others
to show that the superheight is $\geq 2$.
\par\bigskip
\noindent
{\bf Example 3.1.}
Let $K$ be a domain and consider in the domain
$$A=K[X_1,X_2,Y_1,Y_2]/(X_1^kX_2^k+Y_1X_1^{k+1}+Y_2X_2^{k+1})$$
the ideal ${\bf a}=(X_1,X_2)$, $U=D({\bf a})$.
The functions 
$$Z=-Y_2/ X^k_1=(X_2^k+Y_1X_1)/ X^{k+1}_2 \mbox{ and }
W=-Y_1/X^k_2=(X_1^k+Y_2X_2)/ X^{k+1}_1 $$
are definied on $U$ and one has $X_2Z+X_1W=1$, hence $U$ is affine.
\par\smallskip\noindent
This example for $K={\bf Z}$ is the two dimensional case
of the superheight version of the monomial conjecture,
and maybe the easiest way to settle this instance is
by showing the affineness. For another proof see \cite{hochstercan}.
\par\bigskip
\noindent
{\bf Example 3.2.}
Now we look at the prime ideal ${\bf a}=(X_1,X_2)$ in the domain
$$A=K[X_1,X_2,Y_1,Y_2]/(X_1^kX_2^k+Y_1X_1^k+Y_2X_2^{k+1})\, .$$
Consider the morphism $A \longrightarrow A'=K[X_1,X_2]$
given by the substitution $Y_1 \longmapsto -X_2^k,\, Y_2 \longmapsto 0 $.
Then ${\bf a}A'=(X_1,X_2)$ has height two, and $D({\bf a})$
is not affine.
\par
\bigskip
\noindent
{\it Two\--dimensional rings}
\par
\medskip
\noindent
A theorem of Nagata states that
on a normal affine surface
the complement of any (pure one\--dimensional) curve is affine, see
\cite{naghil}. From the proof of this theorem one can get
the following theorem.
\begin{satz}\hspace{-0.5em}{\bf .}
Let $A$ be a two dimensional noetherian ring, ${\bf a} \subseteq A$.
Then $D({\bf a})$ is affine if and only if the noetherian
superheight of $\bf a$ is $\leq 1$.
\end{satz}
{\it Proof}.
Suppose $U=D({\bf a})$ is not affine.
We may assume that $A$ is a two dimensional
noetherian normal and local domain, since the normalization of a 
noetherian two dimensional domain is again noetherian.
$B=\Gamma(U,{\cal O}_X)$ is a Krull domain, and, since
$U$ is not affine, the height of the extended
ideal ${\bf b}={\bf a}B$ is at least two.
By a faithfully flat extension as in \cite{naghil}
one may assume that there exist infinitely many prime elementes in $A$.
Then one can show for a minimal prime ${\bf m}'$ of $\bf b$
that $R=B_{{\bf m}'}$ is the desired two dimensional and noetherian
ring. \hfill $\Box$
\begin{satz}\hspace{-0.5em}{\bf .}
Let $A$ be an excellent two dimensional domain.
The complement of a curve
$Y \subseteq X={\rm Spec}\, A$ is affine if and only
if every component of the preimage of $Y$ in the normalization
$\tilde{X}$ has codimension one. This means that the preimage does not
have isolated points.
\end{satz}
{\it Proof}.
If the preimage $\tilde{Y}$ has pure codimension one,
the theorem of Nagata (which is valid for affine excellent surfaces)
says that $\tilde{Y}$ has an affine complement,
and the theorem of Chevalley says that this holds for $Y$ itself. \hfill $\Box$
\par
\bigskip
\noindent
{\bf Remark.}
Of course, if the normalization is a bijection any complement
of a curve is affine. If this is not the case it is
quite easy to find curves with non\--affine complement. If
$Q,R \in \tilde{X}$ are different points mapping to $P \in X$,
look for curves $Y'$ on $\tilde{X}$
lying generically inside the open set where the normalization
is an open embedding (say $X$ excellent)
and with $Q \in Y',\, R \not \in Y'$.
Then the image $Y$ of $Y'$ cannot have an affine complement,
because the preimage $\tilde{Y}=Y' \cup \{ R\} $ and $R$ is an isolated
point in it.
On $X$ itself we have to look for regular (or at least cuspidal)
curves $C$ through $P$ not totally lying on ${\rm Sing}\, X$.
\par
\bigskip
\noindent
{\it Monoid rings}
\par\medskip
\noindent
Let $M$ be a normal torsion\--free finitely generated monoid with
quotient lattice $\Gamma={\bf Z}M \cong {\bf Z}^d$.
Let $M$ be positive, meaning that $0$ is the only unit of $M$.
Then there exists an embedding with the intersection property,
namely $ M \hookrightarrow {\bf Z}^k$ with
$M= \Gamma \cap {\bf N}^k$, see \cite{bruns}, exc. 6.1.10,
or take the natural embedding given by the
divisor class representation.
Such an embedding yields an inclusion of rings
$$K[M] \hookrightarrow K[{\bf N}^k]=K[T_1,...,T_k] \, , $$
and $K[M]$ is the ring of degree zero
under the $D-$graduation of the polynomial ring given by
${\bf Z}^k \longrightarrow {\bf Z}^k/\Gamma=:D$.
In particular $K[M]$ is a direct summand of $K[T_1,...,T_k]$.
\begin{satz}\hspace{-0.5em}{\bf .}
Let $M$ be a finitely generated torsion free monoid and
$K$ a noetherian factorial domain.
Then there exists a ring extension of finite type
$K[M] \hookrightarrow B$  such that
an open subset $U=D({\bf a}) \subseteq {\rm Spec}\, K[M]$
is affine if and only if ${\rm bight}\, {\bf a}B \leq 1$.
In particular $U$ is affine if and only if
${\rm supht}^{\rm fin} {\bf a} \leq 1$.
\end{satz}
{\it Proof}.
Let $\tilde{M}$ be the normalization of $M$ and
$\tilde{M}={\bf Z}^s \times M'$ with $M'$ positive,
see \cite{bruns}, Theorem 6.1.4 and Prop. 6.1.3. Let
$M' \hookrightarrow {\bf N}^k$ be a representation
with the intersection property. Then the mapping
$$ K[M] \longrightarrow K[\tilde{M}] \cong
K[{\bf Z}^s][M'] \longrightarrow 
K[V_1,...,V_s,V_1^{-1},...,V_s^{-1}][T_1,...,T_k]=B$$
is of finite type.
Let ${\bf a}K[M]$ be an ideal with ${\rm bight}\, {\bf a}B \leq 1$.
Since $B$ is factorial, we know that $D({\bf a}B)$ is affine
and we have to show that this property holds already for
$D({\bf a})$. For a finite extension this is the theorem of
Chevalley, and for a direct summand $A \subseteq B=A \oplus V$
this is true, since $\Gamma (D({\bf a}B),{\cal O}_B)
=\Gamma(D({\bf a}),{\cal O}_A) \oplus \Gamma(D({\bf a}),\tilde{V})$,
and, if $\bf a$ generates the unit ideal in $\Gamma (D({\bf a}B),{\cal O}_B)$,
this is also true in the first component. \hfill $\Box$
\par
\bigskip
\noindent
\section{Affineness, Stein property and superheight one}
\par
\bigskip
\noindent
In the case $K={\bf C}$, we can associate to an algebraic
variety $X$ the corresponding complex space $X^{\rm an}$.
If $X$ is an affine variety, then $X^{\rm an}$ is a Stein space,
see \cite{grauert}, V. \S1 Satz 1.
We will show that the analytic property of being Stein is strong enough
to guarantee that the noetherian superheight is one.
Thus the existence of Stein but non\--affine quasi\--affine
schemes yields directly to non\--affine quasi\--affine varieties
with superheight one.
We consider only separated varieties and complex Hausdorff spaces.
Some results and ideas of this section can also be found in Neeman
(\cite{neeman}) and in Bingener/Storch (\cite{binsto1}).
\par
\bigskip
\noindent
{\bf Definition.}
Let $X$ be a complex space and $Y \subseteq X$ a closed analytic subset.
We define the analytic superheight of $Y$ by
$$ {\rm supht}^{\rm an}\, (Y,X)  =
{\rm sup}\, \{ {\rm codim}_{x'}\,(f^{-1}(Y),X'):
\, x' \in X',\,   f:X' \longrightarrow X \}\,.$$
Here we put ${\rm codim}_x\, (Y,X)={\rm dim}_x \, X-{\rm dim}_x\, Y$
with ${\rm dim}_x\, X={\rm dim}\, {\cal O}_x={\rm dim}\, \hat{{\cal O}}_x$,
see \cite{grauertas}, Kap. II, \S4ff.
If the analytic set $Y_x$ is described
at the point $x \in X$ by the ideal $\bf a$, we have
${\rm codim}_x\, (Y,X)=
{\rm dim}\, ({\cal O}_{X,x}) -{\rm dim}\, ({\cal O}_{X,x}/{\bf a})$.
If $X$ is irreducible this equals the height of the ideal, since
the analytic rings are catenary. 
\begin{lem}\hspace{-0.5em}{\bf .}
Let $Y \subseteq X$ be a closed analytic subset in a 
complex space with $U=X-Y$ Stein.
Then ${\rm supht}^{\rm an}\, (Y,X) \leq 1\, $.
\end{lem}
{\it Proof}.
Let $f:X' \longrightarrow X $ be a morphism of complex spaces
and $x'$ a point of $X'$. Since the codimension is local, we can assume that
$X'$ is Stein.
$f$ factors through the closed graph
$X' \stackrel{ \Gamma_f}{\hookrightarrow} X' \times X 
\stackrel{p_2}{\longrightarrow} X$
and therefore $f^{-1}(X-Y)$ is isomorphic to a closed
subset of $X' \times (X-Y)$. Since $X'$ and $X-Y$ are Stein,
the product $X' \times (X-Y)$
is also Stein and so is $X'-f^{-1}(Y) \subseteq X'$,
see \cite{grauert}, V. \S1 Satz 1.
But the complement of an open Stein subset in a Stein space has
codimension $\leq 1$, see \cite{grauert}, V. \S3, Satz 4. \hfill $\Box$
\begin{satz}\hspace{-0.5em}{\bf .}
Let $X={\rm Spec}\, A$ be an affine algebraic ${\bf C}-$variety and
$V({\bf a})=Y \subseteq X$. Then the algebraic and the analytic
superheight coincide.
$${\rm supht}\, {\bf a}={\rm supht}^{\rm fin}\, {\bf a}
={\rm supht}^{\rm an}\, (Y^{\rm an},X^{\rm an}) \, .$$
\end{satz}
{\it Proof}.
The first equality follows from the theorem of Koh 2.3.
Of course, the analytic superheight is not lower than
the finite algebraic superheight, since we can interpret every
algebraic test variety as an analytic variety, and algebraic and analytic
dimension coincide.
\par\smallskip\noindent
For the converse, let
$f: X' \longrightarrow X^{\rm an}$ be a morphism of complex spaces,
$x' \in X',\, f(x')=x$.
We may suppose that $X'$ is irreducible.
The extended ideal ${\bf a}{\cal O}_{X^{\rm an},x}$
under $A \longrightarrow {\cal O}_{X^{\rm an},x}$
describes the zero set $Y^{\rm an}$ in $x$,
and the preimage $Y'$ in $x'$ is described by ${\bf a}{\cal O}_{X',x'}$.
Since the local rings in a complex space are noetherian, see
\cite{grauertas}, Kap. I, 35.2, Satz 3, and Kap. II, 30, Satz 1,
we have
${\rm codim}_{x'}\, (Y',X')={\rm ht}\, ({\bf a}{\cal O}_{X',x'})
\leq {\rm supht}\, {\bf a}$. \hfill $\Box$
\begin{kor}\hspace{-0.5em}{\bf .}
Let $A$ be a ${\bf C}-$algebra of finite type and
$U=D({\bf a}) \subseteq X$ an open subset with
$U^{\rm an}$ Stein.
Then ${\rm supht}\, {\bf a} \leq 1$.
\end{kor}
Proof.
This follows from the theorem and the lemma. \hfill $\Box$
\begin{kor}\hspace{-0.5em}{\bf .}
Let $A$ be a domain of finite type over ${\bf C}$
and $U \subseteq {\rm Spec}\, A $ an open subset with
$\Gamma(U,{\cal O}_X)$ finitely generated.
Then $U$ is affine if and only if $U$ is Stein.
\end{kor}
{\it Proof}.
The previous corollary shows that the finite superheight is one.
This together with the finiteness of the global ring shows that $U$
is affine.
(For another proof see \cite{binsto1}, 5.1) \hfill $\Box$
\par
\bigskip
\noindent
\section{Superheight one and affineness of two dimen\-sional subschemes}
\par\bigskip
\noindent
Let $U$ be a separated scheme of finite type over a field $K$.
In this section we
study the property that every closed surface in $U$ is affine.
\begin{satz}\hspace{-0.5em}{\bf .}
Let $A$ be a domain of finite type over a field $K$,
$D({\bf a})=U \subseteq X={\rm Spec}\,A$ an open subset.
Then the following are equivalent.
\par\smallskip\noindent
{\rm (1)} ${\rm supht}\, {\bf a} \leq 1$
\par\smallskip\noindent
{\rm (2)} Every closed subvariety of dimension $\leq 2$ of $\, U$ is affine.
\par\smallskip\noindent
If $K={\bf C}$, this is also equivalent to the following.
\par\smallskip\noindent
{\rm (3)} For every closed analytic surface
$S \subseteq X^{\rm an}$ the intersection $S \cap U^{\rm an}$ is Stein.
\end{satz}
{\it Proof}.
Suppose (1) holds. For points and curves the statement (2) is always true,
so let $S \hookrightarrow U$ be a closed reduced irreducible
surface in $U$, and let $S'$ be the closure of $S$ in $X$.
Then $S' \hookrightarrow X$ is again a surface, because the
dimension of an irreducible variety does not change
in passing to a non\-empty open subset.
Let $\tilde{S} \longrightarrow S'$ be the normalization.
The preimage of $Y=V({\bf a})$ under $\tilde{S} \longrightarrow X$ has
due to the superheight property pure codimension one and
hence due to the theorem of Nagata an affine complement.
The theorem of Chevalley shows that the complement
of $S' \cap Y$ is again affine, so
$S=S' \cap U =S'-S' \cap Y$ is affine.
\par\smallskip\noindent
For the converse let
${\rm supht}^{\rm fin}\, {\bf a}={\rm supht}\, {\bf a} \geq 2$.
Then there exists (theorem 2.2) an irreducible
surface ${\rm Spec}\, R=S \subseteq X$ with normalization
$S'={\rm Spec}\, R_{\rm nor}$ such that ${\bf a}R_{\rm nor}$ has big height
2. Thus $D({\bf a}R_{\rm nor})$ and $D({\bf a}R)$ can not be affine,
and $D({\bf a}R)=U \cap S$.
\par\smallskip\noindent
Now let $K={\bf C}$ and suppose (1) holds.
Let $S \subseteq X^{\rm an}$ be a closed analytic surface
with normalization
$f:\tilde{S} \longrightarrow S \hookrightarrow X^{\rm an} $.
Then the codimension of $f^{-1}(Y)$ on the normal surface $\tilde{S}$ is
$\leq 1$, because the algebraic superheight equals the analytic superheight.
The theorem of Simha (this is the analytic analogue to the theorem
of Nagata, see \cite{simha}) says that $\tilde{S} -f^{-1}(Y)$
is a Stein space. This means that the normalization of $U \cap S$ is Stein,
and so due to the analytic version of the theorem of Chevalley
$U \cap S$ itself is Stein.
\par\smallskip\noindent
Now suppose (3) holds and let an algebraic surface $S' \subseteq U$ be given.
We can write $S'=S \cap U$ with a closed algebraic surface $S$ in $X$.
By (3) we know that $S'=S \cap U$ is Stein, so by 4.3 and
3.3 it is affine. \hfill $\Box$
\begin{kor}\hspace{-0.5em}{\bf .}
Let $U$ be a quasi\--affine variety over $K$ such that the ring of global
sections is finitely generated.
If all irreducible closed surfaces of $U$ are affine,
$U$ itself is affine.
\end{kor}
{\it Proof}.
This follows directly from the theorem and theorem 3.2. \hfill $\Box$
\par\bigskip
\noindent
{\bf Remark.}
In case $K={\bf C}$, the last statement of the
theorem is fullfilled if $U$ itself is Stein.
The hypersection (or hypersurface) problem in complex analysis asks the following:
given a Stein space $X$ of dimension $\geq 3$
and an open subset $U \subseteq X$ with the property
that for any analytic hypersurface $S \subseteq X$
the intersection $S \cap U$ is Stein,
is then $U$ itself Stein?
If $U \subseteq X$ is algebraic and ${\rm dim}\, X=3$,
statement (3) of the above theorem
is exactly the condition of the hypersurface problem.
\par\smallskip\noindent
However, the hypersurface problem is now known not to be true in
general, as first shown by an example of Coltoiu and Diederich, see
\cite{diederich2} and \cite{diederich1}.
In section (6) we will give a class of examples of non-Stein
open subsets with superheight one, and 5.1 shows that the assumptions of the
hypersection problem are fulfilled.
\par
\bigskip
\noindent
{\bf Example 5.1.}
The affiness of an open subset cannot be tested (even if the ring
of global sections is finitely generated) with more restrictive classes of
surfaces. The following example shows that homogeneous surfaces do not suffice.
\par\smallskip\noindent
Let $S$ be the projective plane, blown up in one point $P$.
Let $E$ be the exceptional divisor and $C$ a projective line not
passing through the point.
$W=S-(E \cup C)$ is then a punctured affine plane, so $W$ is quasi\--affine
and contains no projective lines.
Let $A$ be a homogenous coordinate ring for $S$, $W=D_+({\bf a})$,
$U=D({\bf a}) \subseteq X={\rm Spec}\, A$.
For an irreducible homogenous surface $V({\bf p})$ in the affine cone $X$
the corresponding projective curve $V_+({\bf p})$ intersects $V_+({\bf a})$,
and therefore $V_+({\bf p}) \cap W$ is affine, hence also
the preimage $V({\bf p}) \cap U $. This means that all homogenous
surfaces inside $U$ are affine.
But not all surfaces in $U$ are affine.
$U$ is the cone over a punctured affine plane
and thus isomorphic to
${\bf A}^\times \times ({\bf A}^2 - \{P\})$.
As a subset of 
${\bf A}^{\times} \times {\bf A}^2$ it has height two,
and this gives a lot of non\--affine surfaces.
\par
\bigskip
\noindent
\section{Non-affine subsets with superheight one}
\par
\bigskip
\noindent
A theorem of Goodman states that on a smooth projective surfaces $S$ an open
subset $U=S-Y$ is affine if and only if there exists an ample effective
divisor $H$ with ${\rm supp}\, H=Y$, \cite{goodman}, \cite{haramp},II.4.2.
A weaker condition on $Y$ and $H$ still implies that $U$ has superheight one.
\begin{satz}\hspace{-0.5em}{\bf .}
Let $S$ be a smooth projective surface over an algebraically closed
field $K$, $Y \subset S$ a curve and $U=S-Y$. Suppose
there exists an effective divisor $H$ with ${\rm supp}\, H=Y$
and with $H.Y_i \geq 0$ for all
irreducible components $Y_i$ of $Y$ and with
$H.C >0$ for all curves $C \not\subseteq Y$.
Then every morphism $T \supseteq T-\{P\} \longrightarrow U$,
where $T$ is a two dimensional normal irreducible affine variety,
is extendible to $T$.
\par\smallskip\noindent
If $S={\rm Proj}\, A$ with a finitely generated graded
$K-$Algebra $A$ and $U=D_+({\bf a})$, then ${\rm supht}\, {\bf a}=1$.
\end{satz}
{\it Proof}.
We have already seen in 2.5 that ${\rm supht}^{\rm fin}\, (U) = 1$
follows from the described extendibility property. Since the cone mapping
is affine, it follows that ${\rm supht}^{\rm fin}\, {\bf a} =1$ and, due to the
theorem of Koh, ${\rm supht}\, {\bf a}=1$.
\par\smallskip\noindent
So let $f:T-\{ P\} \longrightarrow U$ be a morphism of a reduced
irreducible normal affine surface $T$. We may assume that $T -\{P\}$ is
regular.
If $f(T-\{P\})$ is a point, $f$ is of course extendible.
If $f(T - \{ P \}) \subseteq U$ lies inside an irreducible curve
$C \subseteq U$, this curve $C$ is due to the assumption not
projective, hence affine.
Then $f$ corresponds to a ring-homomorphism and is thus
extendible to $C \subseteq U$ as in the proof of theorem 2.6.
So suppose that the image of $f$ is two dimensional and $f$ dominates $U$.
\par\smallskip\noindent
Let $T \hookrightarrow T'$ be an open embedding in a projective
surface with complement $D'$.
Let $p:\tilde{T} \longrightarrow T'$ be a resolution of singularities of $T'$
and a resolution of the undefined points of
$f:T' \supseteq T-\{ P \} \longrightarrow S$, see \cite{haralg}, V.3.8.1
and Theorem V.5.5. So we have an extension
$\tilde{f}: \tilde{T} \longrightarrow S$ of
$f$ on $T-\{ P \} \cong \tilde{T}-p^{-1}(P)-p^{-1}(D')$.
Let $C_1,...,C_n$ be the irreducible one\--dimensional
components of $p^{-1}(P)$ and let
$D_1,...,D_m$ be the components of $p^{-1}(D')$.
\par\smallskip\noindent
Since $\tilde{f}$ is surjective, it induces a mapping
$\tilde{f}^*$ of the divisors (=Cartier\--divisors).
Let
$\tilde{f}^*(H)=C+D$ with $C=k_1C_1+...+k_nC_n$ and $D=l_1D_1+...+l_mD_m$
be the pull\--back of the divisor $H$; there cannot be other components.
$\tilde{f}_*(C)$ is a non\--negative
combination of the $Y_j$, so the assumptions
concerning the intersections of $H$ with its components yield
$$0 \leq \tilde{f}_*(C).H=C.\tilde{f}^*(H)=C.C+C.D=C.C \,\, .$$
But due to \cite{artin}, theorem 2.3, the self intersection number of an
effective divisor $\neq 0$ is negative, if all its
components are (possibly singular) contractible.
Since the components of $C$ are contracted by $p$ to $P$,
we must have $C=0$. So for all components $C_i$ we have
$\tilde{f}(C_i) \not\subseteq Y$.
\par\smallskip\noindent
So the preimage of $Y$ under $\tilde{f}$ contains only some points on
the $C_i$. If $Q \in C_i$ with $\tilde{f}(Q)=R\in Y$, we find
\---since $S$ is regular and hence locally factorial\---
an affine neighbourhood $W$ of $R$ where
$Y$ is described by one function, so the preimage of $Y$ must be a curve,
which is already excluded.
\par\smallskip\noindent
So we conclude that $\tilde{f}(C_i) \subseteq U=S-Y$ for all curves $C_i$.
Since on $U$ there exist no projective curves,
all these curves are contracted by $\tilde{f}$ to a point of $U$, and
to exactly one point, because the $C_i$ are connected.
So $f$ itself is extendible in $P$ as a function to $U$. \hfill $\Box$
\par
\bigskip
\noindent
{\bf Remark.}
We cannot weaken the assumptions in this theorem.
If $Y$ is reducible with $Y^2=0$ and the complement contains no
projective curves, this has no consequence on the
superheight as shown by the example at the end of section
5. If $ Y \cap C = \emptyset $ for some curve $C$, then
$C$ is a projective curve lying inside $U$, and the cone mapping
of this curve is not extendible to the point of the cone.
\par\smallskip\noindent
If $Y$ is irreducible, the condition of the theorem
says $Y^2 \geq 0$ and $S-Y$ contains no projective curve.
In this case we cannot avoid the assumtion $Y^2 \geq 0$.
If $K={\bf C}$ and $Y^2 <0$, one can contract $Y$
onto a (possibly non\--algebraic) complex space, see \cite{grauertmod}.
This contraction yields a mapping back on this complex space defined
outside the contraction point, and this mapping is not extendible.
\begin{kor}\hspace{-0.5em}{\bf .}
Suppose the situation of the theorem holds, but there
exists no ample effective divisor $H$ with support $Y$.
Then $U=D_+({\bf a})$
and the preimage $D({\bf a})$ in the affine cone are not affine, yet
their superheight is one.
This is in particular the case if $\, Y$ is irreducible with $Y^2=0$
or if $\, Y$ is not connected.
\end{kor}
\par
\medskip
\noindent
{\it Proof}.
If $U$ is affine then there exists an ample divisor $H$ with support $Y$,
see \cite{haramp} or \cite{goodman}.
If $U$ is not affine then also the preimage in an affine cone can not
be affine. This can be seen for example by considering the cohomology
of quasi\-coherent sheaves coming from graded modules.
An effective ample divisor has a positive self intersection
number and is connected, see \cite{haralg}, II.6.2. \hfill $\Box$
\par
\bigskip
\noindent
{\bf Remark.}
A problem of Hartshorne (\cite{haralg}, VI.3.4, \cite{umemura}, \cite{vovan})
asks the following: suppose we are given a smooth complete
algebraic surface $S$ over $\bf C$ and an irreducible curve $Y$
intersecting every other curve positively and with self intersection zero.
Is then $S-Y$ Stein?
\par\smallskip\noindent
Our theorem states that in this situation $S-Y$ fulfills all
geometric conditions which would follow from the Stein property,
so at least it is not possible to refute this conjecture by
geometrical means. Furthermore, our theorem relates this problem
to the hypersurface problem: from the assumptions on $Y \subset S$
in the problem of Hartshorne it follows via 6.1 and 5.1 that
the corresponding open subset in an affine cone over $S$ fulfills
the assumptions of the hypersurface problem (the conclusion of both
problems being the same). The original problem of Hartshorne is still open,
Vo Van proves it in \cite{vovan} in the case where $S$ is a ruled surface.
\par
\bigskip
\noindent
We will give some examples of curves on smooth surfaces where the situation
of the corollary (and of the theorem) occurs.
\par\bigskip
\noindent
{\bf Example 6.1.} (see \cite{goodlandman}, 6.10 and \cite{haralg}, V.5.7.3)
Let $K$ be an algebraically closed field,
$Y_0 \subseteq {\bf P}^2_K$ be a smooth curve of degree three, hence an
elliptic curve, and let $P_1,...,P_9$ be nine points on $Y_0$.
Let $S$ be the blown\--up surface of these $9$ points
and let $Y$ be the proper transform of $Y_0$.
The self intersection is $0$.
$Y$ intersects all exceptional divisors and the
intersection with the other curves
is also positive if the points are choosen in such a way that there does
not exist a relation between them in the group structure on $Y_0$.
\par
\bigskip
\noindent
{\bf Example 6.2.} (See \cite{umemura}, \cite{bingener}, \cite{haramp} with the
corrections due to \cite{neeman})
This is the classical example of a non\--affine but Stein surface.
Let
$$0 \longrightarrow {\cal O}_C \longrightarrow {\cal E}
\longrightarrow {\cal O}_C
\longrightarrow 0 $$
be a non split exact sequence of sheaves on an elliptic curve $C$ over $\bf C$,
there $\cal E$ is locally free of rank two.
Let $s:C \longrightarrow S$ be the section in $S={\bf P}({\cal E})$
corresponding to the epimorphism, and put $Y=s(C)$ and $U = S-Y$.
Then $Y$ fulfills the conditions in 6.2, and it is also Stein, the same being
true in the affine cone.
\par
\bigskip
\noindent
We construct a second class of non-affine, quasiaffine schemes with superheight
one. For this, let $R$ be a noetherian normal domain
and let $M$ be a reflexive (finitely generated) $R-$module of rank one,
corresponding to a Weil divisor.
Let $S(M)$ be the symmetric algebra of $M$, and
put $X={\rm Spec}\, S(M)$ with restriction map
$p: X \longrightarrow {\rm Spec}\, R=Y$.
\par\smallskip\noindent
Let $V \subseteq Y$ be an open subset containing the points
of codimension one such that $M$ defines an invertible sheaf $\cal L$ on $V$.
Then $X\mid_V=p^{-1}(V) \longrightarrow V$ is a line bundle.
Its ring of global sections is given by
$$ \Gamma(p^{-1}(V), {\cal O}_X)=
\oplus_{k \geq 0} \Gamma(V,{\cal L}^{k})\,
= \oplus_{k \geq 0} (M^{\otimes k})^{**} \, .$$
If $M={\bf p}$ is a prime ideal of height one,
this ring equals also $A+{\bf p}+{\bf p}^{(2)}+...$.
\par\smallskip\noindent
The zero-section in $X$ defines the closed subscheme $Z=V(S(M)_+)$.
Above $V$ the open subset $U= p^{-1}(V)\cap (X-Z)$ is a $G_m-$fiber bundle,
its ring of global sections is given by
$$\Gamma(U,{\cal O}_X) = \oplus_{k \in {\bf Z}} \Gamma(V,{\cal L}^k) \, .$$
\
A line bundle is trivial if and only if there exists a section without zero
and a $G_m-$fiber bundle is trivial if and only if it has a section.
\begin{satz}
\hspace{-0.5em}{\bf .}
Let $R$ be a noetherian normal domain with a closed
point $P \in {\rm Spec}\, R=Y$ of height $d \geq 2$ such that
$V=Y-\{P \}$ is locally factorial.
Let ${\cal L} \in {\rm Pic}\, V \cong {\rm Cl}\, R$ be a non-torsion element
in ${\rm Cl}\, R_P$.
Let $U$ be the corresponding $G_m-$fiber bundle over $V$.
Then the cohomological dimension of $U$ is $d-1$ and its finite superheight is
$\leq d-1$.
\par\smallskip\noindent
If $P$ is a closed point on a normal affine surface, then $U$ has
finite superheight one, but is not affine.
\end{satz}
{\it Proof}.
For a finitely generated positively graded algebra $S$ over $R$ and a
homogeneous ideal $\bf a$ the cohomological dimension of $D({\bf a})$
and $D_+({\bf a}) \subseteq {\rm Proj}\, S$ are the same.
This follows from the fact that any coherent sheaf on $D_+({\bf a})$
comes from a graded module. We may apply this to $U \longrightarrow V$ and
therefore ${\rm cd}\, U={\rm cd}\, V=d-1$.
\par\smallskip\noindent
Let now $R'$ be a normal noetherian domain of dimension $d$
and let $f: Y'={\rm Spec}\, R' \longrightarrow X={\rm Spec}\, S$
be a morphism of finite type.
We have to show that
$f^{-1}(U) \neq Y'-\{ P' \}$, where $P'$ is a closed point of height
$d$. First observe that $p(f(P'))=P$, for otherwise $p(f(P')) \in W$,
where $W$ is an affine neighbourhood with $X\mid_W$ trivial, and
$f(P') \in Z \cap p^{-1}(W)$, but this is not possible since
${\rm ht}\, P' \geq 2$.
Therefore $g=p \circ f : Y' \longrightarrow Y$ is a morphism of finite
type with $g^{-1}(P)=\{P' \}$, and we have to exclude that
$f: Y'-\{P'\} \longrightarrow X\mid_V $ does not meet $Z$ at all.
But such a mapping would yield a zero-free section
$f' : Y'-\{P'\} \longrightarrow g^*(X\mid_V)$ on the pull back of the line
bundle $X\mid_V$ and this would be trivial, but this is not possible
as the following lemma shows. \hfill $\Box$
\begin{lem}\hspace{-0.5em}{\bf .}
Let $R$ and $R'$ be normal excellent domains
with maximal ideals ${\bf m} $ and ${\bf m}'$ of same height $d \geq 2$.
Let $R \longrightarrow R'$ be a ring homomorphism of finite type with
$V({\bf m}R')=V({\bf m}')$.
Then the kernel of $Cl\, R_P \longrightarrow Cl \, R'_{P'}$
consists of torsion elements.
\end{lem}
{\it Proof}
.
We may assume that $R$ and $R'$ are local, and from $V({\bf m}R')=V({\bf m}')$
we see that also $\hat{R} \longrightarrow \hat{R'}$ is of finite type.
Since we assume excellence, normality is preserved by completion,
and ${\rm Cl}\, R \longrightarrow {\rm Cl}\, \hat{R}$ is injective,
see \cite{fossum}, Cor 6.12.
Thus we may assume that both rings are complete.
Since $R$ and $R'$ have the same dimension and the closed
fiber is zero-dimensional it follows that $R \longrightarrow R'$ 
is quasifinite. Due to (see \cite{EGAII}, 6.2.6) it is already
finite and the result follows by taking the norm. \hfill $\Box$
\par
\bigskip
\noindent
{\bf Example 6.3.}
To construct examples of the desired type
we have to look for affine normal surfaces
$Y={\rm Spec}\,R$ with prime ideals $\bf p$ of height one which are not
torsion at a point $P \in Y$.
One can take for instance the homogeneous
coordinate ring of a smooth projective curve of genus $\geq 1$. If
the curve is elliptic, such divisors are given by points which are not torsion
in the group structure. Another example is given in \cite{bingener}, 2.10. (3).
\par\smallskip\noindent
Examples of such prime ideals were first used by Rees to construct examples
of non-finitely generated rings of global sections.
From the properties established in the theorem it follows by 3.2 that the
global ring of $U$ is not finitely generated.
\par
\bigskip
\noindent
{\bf Remark.}
Take an example as above where
$R$ is a finitely generated normal ${\bf C}-$Algebra of dimension two.
Then $U^{\rm an} \subseteq X^{\rm an}$ is an example of a complement
of a hypersurface in a Stein space,
fulfilling the assumptions in the hypersection problem but not the conclusion.
For in that case it follows from superheight one via 5.1 that
for every closed analytic surface (=hypersurface)
$T \hookrightarrow X^{\rm an}$ the intersection $T \cap U$ is Stein.
However, on a complex manifold $V$ the complement of the zero-section
in a line bundle $L$ can only be Stein in case $V$ itself is Stein, see
\cite{diederich2}, Lemma 3.21.
But here $V=Y-\{ P\}$ is not Stein.
The example of Coltiou/Diederich can be interpreted in this context as
in the context of 6.2 as well.
\par
\bigskip
\noindent
We will discuss a third class of non\--affine schemes with superheight one
arising from tight closure in characterisitic $0$ and related to example 6.2
in another paper.

\end{document}